\documentclass[12pt]{article}
\usepackage{algorithm}
\usepackage{algpseudocode}
\usepackage{amsmath}
\usepackage{amssymb}
\usepackage{bm}
\usepackage{booktabs}
\usepackage{color}
\usepackage{cite}
\usepackage{epsfig}
\usepackage{epstopdf}
\usepackage{flafter}
\usepackage{float}
\usepackage{ifthen}
\usepackage{multirow}
\usepackage{amsfonts}
\usepackage{makecell} 
\usepackage{pifont}
\usepackage{subfigure}
\usepackage{times}
\usepackage{threeparttable}
\usepackage{url} 
\usepackage{subfigure,graphicx}
\usepackage{float}
\usepackage{color,multirow}
\usepackage{makecell}
\usepackage{setspace}
\usepackage{cite} 
\usepackage{url} 
\usepackage{ifthen} 
\usepackage[T1]{fontenc}
\usepackage{graphicx}
\usepackage{changepage}
\usepackage[colorlinks,linkcolor=red,anchorcolor=gray,citecolor=green,urlcolor=blue]{hyperref}
\usepackage{graphics}
\usepackage{booktabs}

\usepackage{amsfonts}
\usepackage{color}
\usepackage{amssymb}
\usepackage{amsmath}
\usepackage{indentfirst}

\def\sqr#1#2{{\vcenter{\vbox{\hrule height.#2pt
				\hbox{\vrule width.#2pt height#1pt \kern#1pt \vrule width.#2pt}
				\hrule height.#2pt}}}}
\def\signed #1{{\unskip\nobreak\hfil\penalty50
		\hskip2em\hbox{}\nobreak\hfil#1
		\parfillskip=0pt \finalhyphendemerits=0 \par}}
\def\endpf{\signed {$\sqr69$}}

\def\dbR{{\mathop{\rm l\negthinspace R}}}

\def\3n{\negthinspace \negthinspace \negthinspace }
\def\2n{\negthinspace \negthinspace }
\def\1n{\negthinspace }

\def\dbE{\mathbb{E}}
\def\dbF{{\mathbb{F}}}
\def\BF{{\bf F}}

\def\ds{\displaystyle}

\def\dbN{{\mathop{\rm l\negthinspace N}}}
\def\dbP{{\mathbb{P}}}
\def\dbR{{\mathbb{R}}}
\def\mE{{\mathbb{E}}}
\def\={\buildrel \triangle \over =}

\def\a{\alpha}
\def\b{\beta}

\def\d{\delta}
\def\e{\varepsilon}

\def\k{\kappa}
\def\l{\lambda}
\def\m{\mu}
\def\n{\nabla}
\def\si{\sigma}
\def\t{\times}

\def\th{\theta}
\def\om{\omega}

\def\i{\infty}
\def\ns{\noalign{\ss} }

\def\G{\Gamma}
\def\D{\Delta}
\def\Th{\Theta}

\def\Si{\Sigma}

\def\Om{\Omega}

\def\cA{{\cal A}}
\def\cB{{\cal B}}
\def\cC{{\cal C}}

\def\cF{{\cal F}}

\def\cH{{\cal H}}

\def\cJ{{\cal J}}
\def\cK{{\cal K}}

\def\cM{{\cal M}}
\def\cN{{\cal N}}

\def\no{\noindent}

\def\ss{\smallskip}
\def\ms{\medskip}
\def\bs{\bigskip}
\def\q{\quad}
\def\qq{\qquad}

\def\max{\mathop{\rm max}}
\def\min{\mathop{\rm min}}
\def\exp{\mathop{\rm exp}}

\def\pa{\partial}

\def\cd{\cdot}

\def\div{\hbox{\rm div$\,$}}

\def\|{\Big |}
\def\({\Big (}
\def\){\Big )}
\def\[{\Big[}
\def\]{\Big]}
\def\be{\begin{equation*}}
	\def\bel{\begin{equation}\label}
		\def\ee{\end{equation}}
	\def\bt{\begin{theorem}}
		\def\bcd{\begin{condition}}
			\def\ecd{\end{condition}}
		\def\et{\end{theorem}}
	\def\bc{\begin{corollary}}
		\def\ec{\end{corollary}}
	\def\bde{\begin{definition}}
		\def\ede{\end{definition}}
	\def\bl{\begin{lemma}}
		\def\el{\end{lemma}}
	\def\bp{\begin{proposition}}
		\def\ep{\end{proposition}}
	\def\br{\begin{remark}}
		\def\er{\end{remark}}
	\def\ba{\begin{array}}
		\def\ea{\end{array}}
	\def\ed{\end{document}}
\def\ns{\noalign{\ms}}
\def\ds{\displaystyle}

\def\square#1{\vbox{\hrule\hbox{\vrule height#1%
			\kern#1\vrule}\hrule}}
\def\rectangle#1#2{\vbox{\hrule\hbox{\vrule height#1%
			\kern#2\vrule}\hrule}}

\font\tenbb=msbm10 \font\sevenbb=msbm7 \font\fivebb=msbm5

\newfam\bbfam
\scriptscriptfont\bbfam=\fivebb \textfont\bbfam=\tenbb
\scriptfont\bbfam=\sevenbb

\newtheorem{lemma}{Lemma}[section]
\newtheorem{remark}{Remark}[section]

\newtheorem{theorem}{Theorem}[section]
\newtheorem{corollary}{Corollary}[section]

\newtheorem{definition}{Definition}[section]
\newtheorem{proposition}{Proposition}[section]
\newtheorem{condition}{Condition}[section]

\allowdisplaybreaks

\makeatletter

\@addtoreset{equation}{section}
\makeatother
\begin{document}
	\title{Exact Controllability for a Refined Stochastic Wave Equation \thanks{This work is partially supported by the NSF of China
			under grants 11971333, 11931011  and 12025105, and by the Science Development Project of Sichuan University under
			grant 2020SCUNL201. }}
	\author{  Zhonghua Liao\thanks{School of Mathematics,  Sichuan University,
			Chengdu 610064,  China. E-mail address: zhonghualiao@yeah.net} \  and\  Qi L\"u\thanks{School
			of Mathematics, Sichuan University, Chengdu
			610064, China. E-mail address: lu@scu.edu.cn.}}
	\date{}
	\maketitle
	\begin{abstract}
		In this paper, we obtain the exact controllability for a refined stochastic wave equation with three controls by  establishing a novel Carleman estimate for a backward hyperbolic-like operator. Compared with the known result (\cite[Chapter 10]{QLXZ}), the novelty of this paper is twofold: (1) Our model contains the effects in the drift terms when we put controls directly in the diffusion terms,  which is more sensible for practical applications; (2) We provide an explicit description of the waiting time which is sharp in the case of dimension one and is independent of the coefficients of lower terms.
	\end{abstract}
	
	\no{\bf 2020 Mathematics Subject Classification}. 93B05, 93B07 
	
	\bs
	
	\no{\bf Key Words}. stochastic wave equation, exact controllability, Carleman estimate.
	
	\section{Introduction}

	\par Let $T>0$, and $G\subset \dbR^n$ ($n\in\dbN $) a  bounded  domain with a $C^2$ boundary $\G$. Let $\G_0$ be a nonempty subset of $\G$, which will be given later.  Write
	\bel{n}
	\ba{ll}
	\ds Q=(0,T)\t G, \q \Si=(0,T)\t \G, \q \Si_0=(0,T)\t \G_0.
	\ea 
	\ee 
	Let $(\Om, \cF, \BF,\dbP)(\mbox{with }\BF\=\{\cF_t\}_{t\in [0,T]}$ ) be a  filtered probability space, on which a one dimensional standard Brownian motion $W(\cd)$ is defined, and $\BF $ is the natural filtration generated by $W(\cd)$. Denote by $\dbF$ the progressive $\si$-field of $\BF$.

	Let $H$ be a Banach space. Denote by
	$L^{2}_{\dbF}(0,T;H)$ the Banach space
	consisting of all $H$-valued and $\dbF$-adapted
	processes $X(\cdot)$ such that
	$\mathbb{E}(|X(\cdot)|^2_{L^2(0,T;H)}) <
	\infty$; by $L^{\infty}_{\dbF}(0,T;H)$ the
	Banach space consisting of all $H$-valued and
	$\dbF$-adapted, essentially bounded processes;
	and by $C_{\dbF}([0,T];L^p(\Omega;H))$
	the Banach space consisting of all $H$-valued
	and $\dbF$-adapted processes $X(\cdot)$ such
	that $X(\cd):[0,T] \to L^p_{\cF_T}(\Omega;H)$ is
	continuous  ($p\in[1,\infty]$). Similarly, one
	can define $C^{k}_{\dbF}([0,T];L^{2}(\Om;H))$
	for any positive integer $k$. All of these
	spaces are endowed with their canonical norms.
	
	\par Consider the following  controlled refined stochastic wave equation with three controls:
	\bel{system1}
	\left\{\ba{ll}
	\ds dy=(\hat y+a_5 f) dt+(a_3 y+f)d W(t) & \mbox{ in }Q,\\
	\ns\ds d\hat y-\D ydt=(a_1 y+a_4 g)dt  +(a_2 y+g)dW(t)&\mbox{ in }Q,\\
	\ns\ds y=\chi_{\G_0}h&\mbox{ on }\Si,\\
	\ns\ds y(0)=y_0,~\hat y(0)=\hat y_0&\mbox{ in }G.
	\ea \right.
	\ee 
	Here   $(y_0,\hat y_0)\in L^2(G)\t H^{-1}(G)$, $(y,\hat y)$ is the state variable, $a_i \in L_\dbF^\i (0,T;L^\i (G))$ ($i=1,2,3,5$) and $ a_4\in L_\dbF^\i (0,T; W_0^{1,\i} (G))$, and $f\in L^2_\dbF(0,T;L^2(G))$, $g\in L_\dbF^2(0,T;H^{-1}(G))$ and $h\in L_\dbF^2(0,T;L^2(\G))$ are three controls. 
	
	{{	
			\begin{remark}
				If $f=a_3\equiv 0$, then the system (\ref{system1}) is reduced to the classical stochastic wave equation  (e.g. \cite{minicourse}). However, the exact controllability of that system fails for any $T>0$ and $\cH \subset (0,T)\t \Si$ even if  the control $g$ acts on everywhere (\cite{QLXZ1}). Motivated by this feature, L\"u and Zhang (\cite{QLXZ1}) present a refined stochastic wave equation ((\ref{system1}) without controls) under a "Stochastic Newton's Law" in \cite{EN} and study the exact controllability problem under the condition that $a_5=0$.   Both the classical stochastic wave equation and the refined
				one are mathematical models for the vibration of strings and
				membranes perturbed by random forces, as well as the propagation of
				waves in random environment. According to the
				main result in this paper, the controllability property of the refined one is
				better than the classical one. Consequently, it may serve as a
				better model when people study control problems for some stochastic systems.
			\end{remark}
			\begin{remark}	
				We would like to point out  that, both $a_5fdt$  and $ a_4gdt$ reflect  effects on the drift terms, when controls $f$ and $g$ are added to the diffusion terms. They are side effects of the controls $f$ and
				$g$ we put in the diffusion rather than the controls we want to put
				on the system. The reason for the appearance of them is that, generally
				speaking, if one put a control in the diffusion term, it will
				influence the drift term indirectly. Here we assume the effect is
				linear. Hence, the terms $a_5fdt$ and $ a_4gdt$ appear.
				
				Usually,  the side effects is an `` uncertain part'' of the system, i.e., one cannot choose $a_4$
				and $a_5$. In such case, when we design the
				controls $f$ and $g$, we should try to eliminate
				the influence of $a_5fdt$ and $ a_4gdt$. This
				makes the controllability problem be more
				complicate than the case $a_4=a_5=0$. In [15], the authors assume that $ a_5\equiv 0 $ to simplify the problem. We improve the Carleman estimate for backward stochastic wave equation to handle the general case $ a_5\neq 0 $, which is one of the  main contribution of this paper.
			\end{remark}
	}}

	The control system  (\ref{system1}) is a stochastic PDE with a nonhomogeneous  boundary condition. Its solution is understood in the sense of transposition (see \cite[Section 10.2]{QLXZ} for example). For the readers' convenience, let us first recall the definition of the solution to (\ref{system1}).
	To this end, we consider
	the following backward stochastic wave equation:
	\begin{equation}\label{bsystem1}
		\left\{
		\begin{array}{ll}
			\ds dz=\hat zdt + Z dW(t) &\mbox{ in } Q_\tau,\\
			\ns\ds d\hat z - \D z dt = (a_1 z +a_2 Z-a_3\widehat Z)dt  +  \widehat Z  dW(t)  &\mbox{ in } Q_\tau,\\
			\ns\ds z = 0 &\mbox{ on }  \Si_\tau,\\
			\ns\ds z(\tau) = z^{\tau},\q \hat z(\tau) = \hat
			z^{\tau} &\mbox{ in } G,
		\end{array}
		\right.
	\end{equation}
	where $\tau\in (0,T]$, $Q_\tau\=(0,\tau)\times
	G$, $\Si_\tau\=(0,\tau)\times \G$  and  $(z^{\tau},
	\hat z^{\tau}) \in L^2_{\cF_\tau}(\Om;H_0^1(G)\t
	L^2(G))$. 
	
	\ss
	
	We first
	recall the definition of the solution to
	\eqref{bsystem1}.
	\begin{definition}\label{def bt sol}
		A quadruple of stochastic processes $(z,Z,\hat
		z,\widehat Z)\in  C_{\dbF}([0,\tau];L^2(\Omega,H^1_0(G))) \times
		L^2_{\dbF}(0,\tau;H^1_0(G))\times 
		C_{\dbF}([0,\tau];L^2(\Omega,L^2(G))) \times
		L^2_{\dbF}(0,\tau;L^2(G))$ is called a solution
		of the system \eqref{bsystem1} if for every
		$\psi\in C_0^\infty(G)$ and a.e. $(t,\om)\in
		[0,\tau]\times\Om$, it holds that
		\begin{equation}\label{def id2}
			z^\tau(x)-z(t,x)=\int_t^\tau\hat z(s,x)ds +
			\int_t^\tau Z(s,x) dW(s)
		\end{equation}
		and
		\begin{equation}\label{def id2-1}
			\begin{array}{ll}
				\ds \int_G  \hat z^{\tau}(x)\psi(x)dx - \int_G
				\hat z(t,x)\psi(x)dx  + \int_t^\tau\int_G
				\n z(s,x)\cd \n 
				\psi(x)dxds
				\\ \ns\ds = \int_t^\tau \int_G
				\big(a_1(s,x) z(s,x)
				+ a_2(s,x)Z(s,x)- a_3(s,x)\widehat
				Z(s,x)\big)\psi(x)dxds\\
				\ns\ds\q + \int_t^\tau\int_G \widehat Z(s,x)
				\psi(x)dxdW(s).
			\end{array}
		\end{equation}
	\end{definition}
	
	\medskip
	
	We have the the following known well-posedness
	result for \eqref{bsystem1} (e.g.
	\cite[Section 4.2]{QLXZ}).

	\begin{lemma}\label{well posed1}
		For any $(z^{\tau},\hat z^{\tau})\in
		L^2_{\cF_\tau}(\Om;H_0^1(G))\times
		L^2_{\cF_\tau}(\Om;L^2(G))$, the system
		\eqref{bsystem1} admits a unique solution $
		(z,Z,\hat z,\widehat Z)$.
		Moreover,
		\begin{equation}\label{best1}
			\begin{array}{ll}\ds
				|z|_{C_\dbF([0,\tau];L^2(\Omega,H^1_0(G)))} +
				|Z|_{L^2_\dbF(0,\tau;H^1_0(G))}+|\hat
				z|_{C_\dbF([0,\tau];L^2(\Omega,L^2(G)))} + |\widehat
				Z|_{L^2_\dbF(0,\tau;L^2(G))}\\
				\ns\ds  \leq Ce^{Cr_1}
				\big(|z^{\tau}|_{L^2_{\cF_\tau}(\Om;H_0^1(G))}+|\hat
				z^{\tau}|_{L^2_{\cF_\tau}(\Om;L^2(G))}\big),
			\end{array}
		\end{equation}
		where\vspace{-0.2cm}
		$$ r_1\= \sum_{i=1}^3|a_i|^2_{L^\infty_{\dbF}(0,T;L^{\infty}(G))},
		$$
		and $ C $ is a constant independent of $ r_1 $ and  $ (z,Z,\hat z,\widehat Z) $. 
	\end{lemma}

	Further,   	solutions to \eqref{bsystem1} satisfies the following hidden regularity.

	\begin{lemma}\label{prop-hid}\cite[Proposition 10.6]{QLXZ}
		Let $(z^{\tau},\hat z^{\tau})\in
		L^2_{\cF_\tau}(\Om;H_0^1(G))\times
		L^2_{\cF_\tau}(\Om;L^2(G))$. Then the solution
		$(z,Z,\hat z,$ $\widehat Z)$ of \eqref{bsystem1}
		  satisfies $\frac{\pa z}{\pa\nu}\big|_{\G}\in
		L^2_{\dbF}(0,\tau;L^2(\G))$. Furthermore,
		\begin{equation}\label{hid-eq1}
			\Big| \frac{\pa z}{\pa\nu}
			\Big|_{L^2_\dbF(0,\tau;L^2(\G))} \leq
			Ce^{Cr_1}\big(|z^{\tau}|_{L^2_{\cF_\tau}(\Om;H_0^1(G))}+|\hat
			z^{\tau}|_{L^2_{\cF_\tau}(\Om;L^2(G))}\big),
		\end{equation}
		where the constant $C$ is independent of $\tau$.
	\end{lemma}

	Now we can give the definition of the transposition solution to \eqref{system1}.
	\begin{definition}\label{1-def1}
		A pair of stochastic processes
		$ (y,\hat y) \in
		C_{\dbF}([0,T];L^2(\Om;L^2(G)))\,\times\,
		C_{\dbF}([0,T];$ $L^2(\Om;$ $H^{-1}(G))) $
		is a transposition solution to \eqref{system1}
		if for any $\tau\in (0,T]$ and $(z^{\tau},\hat
		z^{\tau})\in L^2_{\cF_\tau}(\Om;$ $H_0^1(G))\times
		L^2_{\cF_\tau}(\Om;$ $L^2(G))$, it holds that
		\begin{equation}\label{def id1}
			\begin{array}{ll}
				\ds \mE \langle \hat
				y(\tau),z^{\tau}\rangle_{H^{-1}(G),H^1_0(G)} - 
				\langle \hat
				y_0,z(0)\rangle_{H^{-1}(G),H^1_0(G)} \\
				\ns\ds\q -
				\mE\langle y(\tau),\hat z^{\tau}
				\rangle_{L^2(G)} + \langle y_0,\hat
				z(0)\rangle_{L^2(G)}
				\\ \ns\ds = -
				\mE\int_0^\tau\langle f,a_5\hat z+\widehat Z
				\rangle_{L^2(G)} dt + \mE\int_0^\tau \langle
				g,a_4z+Z\rangle_{H^{-1}(G),H_0^1(G)}
				dt  \\
				\ns\ds \q-\mE\int_0^\tau \int_{\G_0}\frac{\pa
					z}{\pa\nu}hd\G ds .
			\end{array}
		\end{equation}
		Here $(z,Z,\hat z,\widehat Z)$ solves
		\eqref{bsystem1} with $(z^{\tau},\hat
		z^{\tau})\in L^2_{\cF_\tau}(\Om; H_0^1(G))\times
		L^2_{\cF_\tau}(\Om; L^2(G))$.
		
	\end{definition}
	\begin{remark}
		For readers who are not familiar with the notion of transpostion soluiton, we point out the following fact:	
		
		When $h=0$, the control system  \eqref{system1} is a homogeneous boundary value
		problem. By the classical theory for stochastic
		evolution equations, \eqref{system1}  admit a unique weak
		solution 
		$(y,\hat y)\!\in\! C_\dbF([0,T];L^2(\Om; 
		L^2(G)))\times
		C_\dbF([0,T];L^2(\Om;H^{-1}(G)))$ (e.g. \cite[Chapter 3]{QLXZ}) . It follows
		from It\^o's formula that these solutions are
		respectively transposition solution  to
		\eqref{system1}. Then, by
		the uniqueness of the transposition solution to
		\eqref{system1}, we know
		that it is also
		the weak solution to \eqref{system1}. In such sense, the transposition solution to \eqref{system1} is a natural generalization of the weak solution  to \eqref{system1}. 
	\end{remark}
	
	With the aid of Lemma \ref{prop-hid}, by the well-posedness result for stochastic evolution equations with unbounded control operators (e.g.,\cite[Theorem 7.12]{QLXZ}), we have the following well-posedness result for the control system \eqref{system1}.

	\begin{proposition}\label{well posed1-1}
		For each $(y_0,\hat y_0)\in L^2(G)\times
		H^{-1}(G)$, the system \eqref{system1} admits a
		unique transposition solution $(y,\hat y)$.
		Moreover,
		\begin{equation}\label{well posed est1}
			\begin{array}{ll}\ds
				|(y,\hat y)|_{C_{\dbF}([0,T];L^2(\Om;L^2(G)))\times C_{\dbF}([0,T];L^2(\Om;H^{-1}(G)))}\\
				\ns\ds \leq Ce^{Cr_2}\big( |y_0|_{L^2(G)} +
				|\hat y_0|_{H^{-1}(G)} +
				|f|_{L^2_{\dbF}(0,T;L^2(G))} +
				|g|_{L^2_{\dbF}(0,T;H^{-1}(G))}+
				|h|_{L^2_{\dbF}(0,T;L^2(\G_{0}))}\big),
			\end{array}
		\end{equation}
		where
		\begin{equation}\label{r2}
			r_2 \ =
			\sum_{k=1}^3|a_k|^2_{L_{\dbF}^{\infty}(0,T;L^{\infty}(G))}+|a_5|^2_{L_{\dbF}^{\infty}(0,T;L^{\infty}(G))}+|a_4|^2_{L_{\dbF}^{\infty}(0,T;W_0^{1,\infty}(G))}.
		\end{equation}
	\end{proposition}

	Now we can give the definition of the exact controllability for (\ref{system1}).
	\begin{definition}
		The system (\ref{system1}) is called exactly controllable at time $T$ if for any $(y_0,\hat y_0)\in L^2(G)\t H^{-1}(G)$ and $(y_1,\hat y_1)\in L_{\cF_T}^2(\Om;L^2(G))\t L^2_{\cF_T}(\Om;H^{-1}(G))$, one can find  controls $(f,g,h)$ $\in L^2_\dbF(0,T;L^2(G))\t L^2_\dbF(0,T;H^{-1}(G))\t L^2_\dbF (0,T; L^2(\G_0))$ such that the corresponding solution $(y,\hat y)$ to (\ref{system1}) satisfies that $(y(T),\hat y(T))=(y_1,\hat y_1)$.
	\end{definition}
	%
	

	Controllability  problems of
	deterministic partial differential equations
	(PDEs for short) are extensively studied for a long time.  In particular, we refer the readers to \cite{AI,Bardos-Lebeau-Rauch1,
		JLL, DLR, Zhangxu4, XZ, EZ} and the rich
	reference therein for exact
	controllability problems for deterministic wave
	equations. Compared with the determinsitic counterpart, controllability problems for
	stochastic PDEs are much less investigated. Fortunately, those problems   attract more and more attention in recent years. Some interesting works have been done for controllability problems of
	stochastic PDEs. In this respect, we refer to
	\cite{barbu1, FL1, GCL, LiuY, Lu0, Luqi8, Lu2022b,Lu2022,TZ,Yan2022,Yu2022} for some
	results on the controllability of stochastic
	parabolic, complex Ginzburg-Landau, degenerate heat,
	Kuramoto-Sivashinsky, Schr\"odinger,
	transport and beam equations. 
	
	To the best of  our
	knowledge, \cite{QLXZ1} (see also \cite[Chapter 10]{QLXZ}) is the only work 
	concerning the exact controllability of
	stochastic wave equation. In that paper, the authors assume that $a_5=0$, that is, the control in the diffusion term of the first equation in \eqref{system1} does not effect the drift term. This assumption is restrictive. Furthermore, the time $T$ in \cite{QLXZ} for the exact controllability of \eqref{system1} depends on $a_4$. This is not natural for the controllability of both deterministic and stochastic PDEs. In this paper, we prove the exact controllability of  \eqref{system1} and drop the assumption that $a_5=0$ and the dependence of the control time $T$ on $a_4$. 
	
	Generally speaking, one can find the following four main methods for
	solving the exact controllability problem of deterministic wave
	equations:
	\begin{itemize}
		\item The first one is  the Ingham type
		inequality method (e.g., \cite{AI}), which works well
		for wave equations involved in some special
		domains, i.e., intervals and rectangles.
		However,  it is very hard to be
		applied to equations in general domains.
		
		\item The second one is the 
		multiplier method (e.g., \cite{JLL}), which is used
		to treat wave equations with time independent
		lower order terms. It seems that it
		does not work for our problem since the
		coefficients of lower order terms  are  time
		dependent.
		
		\item The third one is the microlocal analysis
		approach (e.g., \cite{Bardos-Lebeau-Rauch1}).  There may be lots of obstacles needed
		to be surmounted if one wants to utilize this
		approach to study stochastic control problems due to the lack of theory for propagation of singularities of stochastic wave equations.

		\item The last one is the global Carleman estimate
		(e.g., \cite{Zhangxu4}), which can be used  to handle wave equations with general  lower order terms and to give explicit bounds on the
		observability constant/control cost in terms of
		the potentials entering in it.
	\end{itemize}
	{{On the other hand, as far as we know, there are mainly two methods to study the controllability problem for stochastic PDEs (more details can be found in \cite{QLXZ}): 
			\begin{itemize}
				\item  The first one is Lebeau-Robbiano iteration method (e.g., \cite{Lu0}), which can be used to solve null controllability problem for stochastic parabolic equations. This method relies on the fast decay of the high frequency component of the solution to the control system. As a result, it  cannot be applied to handle controllability problem for stochastic wave equations.  
				
				\item  The second one is the global Carleman estimate. Recently, there are many works employing this method to study controllability problems of stochastic PDEs (e.g., \cite{barbu1, FL1, GCL, LiuY, Luqi8, Lu2022b,Lu2022,TZ,Yan2022,Yu2022}).  
			\end{itemize}
			In this paper, we use the Carleman estimate to get the desired exact contrllability result, i.e., we establish a novel Carleman estimate for backward hyperbolic-like operator, in which the effects on the drift terms can be seen as lower terms. Compared with \cite{QLXZ1}, we improve the Carleman estimate. As a result, we can deal with the case   $a_5\neq 0$. Moreover, by an explicitly choosing of cut-off function in our proof, the time $T$ can be independent on $a_4$ and $a_5$, which means we can still expect a specific time $T$ to get the exact controllability even though the effects on the drift terms are unknown.}}
	
	\par The rest of this paper is organized as follows. We state our main result in Section \ref{sec-main}. Then,  we  establish a explicit Carleman estimate for backward stochastic wave equations in Section \ref{sec-car}. At last, we prove the observability estimate in Section \ref{sec-ob}.

	\section{Main result}\label{sec-main}

	\par We first give the choice of the time $T$ for the exact controllability of \eqref{system1}.  
	Let $x_0=(x_{10},x_{20},\cdots,$ $x_{n0})\in \dbR^n\backslash \overline G $ be a fixed point satisfying that for any $x=(x_1,x_2,\cdots, x_n)$, 
	\bel{za1}
	\min_{x\in \overline G, i=1,2,\cdots,n}|x_i-x_{i0}|>0.
	\ee 
	Set
	\bel{d10}
	\ba{ll}
	\ds \G_0\=\Big\{ x\in \G\|~ (x-x_0)\cd\nu>0\Big\},\q  \Si_0\=(0,T)\t \G_0,
	\ea
	\ee 
	where $\nu $ denotes the unit outward normal vector  of $G$. Put
	\begin{equation}\label{tr1}
		\ba{ll}
		\ds R_1= \max_{x\in \overline G, i=1,2,\cdots, n}|x_i-x_{0i}|,\\
		\ns\ds T^*=2\sqrt{n}R_1 \max_{x\in \overline G}\sqrt \frac{\max_{i=1,2,\cdots,n}(x_i-x_{0i})^2}{\min_{i=1,2,\cdots,n}(x_i-x_{0i})^2}.
		\ea 
	\end{equation}

	We have the following exact controllability result.
	\begin{theorem}\label{th1}
		For any $T>T^*$, the control system (\ref{system1}) is exactly controllable at $T$ with $\cH=\Si_0$.
	\end{theorem}
	\begin{remark} 
		Due to the finite speed of propagation for solutions to (\ref{system1}), we know that $T$ should be large enough to guarantee the observability estimate. A sharp time $T$ for deterministic wave equations is given in \cite{Bardos-Lebeau-Rauch1}. When $n=1$, we have $T^*=2R_1$, which coincides the time in \cite{Bardos-Lebeau-Rauch1}. However, when $n>1$, we believe $T^*$ given in \eqref{tr1} is not sharp. How to improve $T^*$ is an interesting problem. 
	\end{remark}
	\begin{remark} 
		We employ three controls to achieve the exact controllability of \eqref{system1}. It seems that we use too many controls. However, according to \cite[Theorem 10.9]{QLXZ}, we know that all these controls are necessary.
	\end{remark}

	To prove Theorem \ref{th1}, following the standard duality argument (e.g., \cite[Section 7.3]{QLXZ}), we only need to prove the following observability estimate:
	\begin{theorem}\label{1.1}
		For any $T>T^*$, there exists a constant $C>0$ such that for any solution to \eqref{bsystem1} with $\tau=T$  
		and $(z^T,\hat z^T)\in L_{\cF_T}^2(\Om;H_0^1(G)\t L^2(G))$, it holds that
		\bel{3}
		\ba{ll}
		\ds |(z^T,\hat z^T)|^2_{L^2_{\cF_T}(\Om;H_0^1(G)\t L^2(G))}\\
		\ns\ds \le C \exp(\exp(\exp(e^{Cr_2})))\(\dbE\int_{\Si_0} \|\frac{\pa z}{\pa \nu}\|^2d\G dt +|a_4z+Z|^2_{L^2_{\dbF}(0,T;H_0^1(G))} +|a_5\hat z+\widehat Z|^2_{L^2_\dbF(0,T;L^2(G))}\), 
		\ea 
		\ee 
		where $r_2 $ is given in (\ref{r2}). 
	\end{theorem}

	As said before, in this paper, we prove Theorem \ref{1.1} by a global Carleman estimate. To this end, we borrow some idea from the
	proof of \cite[Theorem 10.19]{QLXZ}. However,  
	we cannot simply mimic the method in
	\cite{QLXZ}. If we do so, the terms concerning $a_5$ cannot be handled since it will have the same order as the energy term $\int_Q|\nabla z|^2 dxdt$. Further, if we follow the method in \cite{QLXZ} directly, the time $T$ will depend on the coefficients of the lower order terms, which is unnatural for the control of PDEs and  stochastic PDEs and may lead the time $T$ for the exact controllability of \eqref{system1} to be very large.  

	\begin{remark}
		Theorem \ref{th1} concludes the exact controllability
		of the system \eqref{system1} with a triple
		$(f,g,h)$, where $g\in L^2_\dbF(0,T;H^{-1}(G))$. 
		Such control is very irregular.  It is very interesting
		to see whether \eqref{system1} is exactly
		controllable when $g\in L^2_\dbF(0,T;L^2(G))$.
		By duality argument, one can show that this is
		equivalent to the following observability
		estimate:
		\begin{equation}\label{4.1-eq1.1}
			\begin{array}{ll}\ds
				|(z^T,\hat
				z^T)|_{L^2_{\cF_T}(\Om;H_0^1(G))\times
					L^2_{\cF_T}(\Om;L^2(G))}\\
				\ns\ds\leq C \exp(\exp(\exp(e^{Cr_2})))\(\dbE\int_{\Si_0} \|\frac{\pa z}{\pa \nu}\|^2d\G dt +|a_4z+Z|^2_{L^2_\dbF(0,T;L^2(G))} +|a_5\hat z+\widehat Z|^2_{L^2_\dbF(0,T;L^2(G))}\),
			\end{array}
		\end{equation}
		where $(z,Z,\hat z,\widehat Z)$ is the solution
		to \eqref{bsystem1} with $\tau=T$ and the final
		datum $(z^T,\hat z^T)$. However, we do
		not know how to prove \eqref{4.1-eq1.1}.
	\end{remark}

	\section{Carleman estimate for backward stochastic wave equation}\label{sec-car}
	
	In this section, we establish a Carleman estimate for a backward stochastic wave equation. To this end, we first recall the following fundamental identity for stochastic hyperbolic-like operators, which is a special case of  \cite[Lemma 6.1]{QLXZ1}.

	\begin{lemma}\label{lem}
		Let $u$ be  an  $H^2(\dbR^n)$ valued It\^o process  and $\hat u$ be an $L^2(\dbR^n)$ valued It\^o process such that 
		\bel{zz}
		d u =\hat u dt +UdW(t) \qq \mbox{ in }(0,T)\t \dbR^n
		\ee
		for some $U\in L^2_\dbF(0,T;H^1(\dbR^n))$. Set $\th=e^\ell$, $v=\th u$ and $\hat v=\th \hat u+\ell_t v$. Then, for $a.e.$ $x\in \dbR^n$ and $\dbP$-a.s,
		\begin{eqnarray}\label{pointwise_identity}
			&&\th\big(-2\ell_t\hat v+2\n\ell\cd \n v+\Psi v\big)\big(d\hat u-\D u dt\big)+\div V dt+d \cM \nonumber\\
			&& =\[ \big(\ell_{tt}+\D \ell -\Psi\big)\hat v^2+\big(\ell_{tt}-\D\ell+\Psi\big)|\n v|^2+2\sum_{j,k=1}^n\ell_{x_jx_k}v_{x_j}v_{x_k}\\
			&& \q  -4\n\ell_t\cd\n v\hat v+\cB v^2+\big(-2\ell_t\hat v+2\n \ell\cd\n v+\Psi v\big)^2\]dt+\cN\nonumber\\
			&&\q -\Big\{\th\big(-2\ell_t\hat v+2\n\ell\cd\n v+\Psi v\big)\ell_t U +\big[2\n(\th U)\cd\n \ell\hat v-\th\Psi_tv U+\th\Psi\hat v U\big] \nonumber\\
			&&\q -2\big[ \n v\cd\n(\th U)+\th \cA v U\big]\ell_t\Big\}dW(t), \nonumber
		\end{eqnarray} 
		where
		\bel{VMBA}
		\left\{\ba{ll}
		\ds \cA \= \ell_t^2-\ell_{tt}-|\n\ell|^2+\D \ell-\Psi,\\
		\ns\ds \cM\=\ell_t|\n v|^2+\ell_t \hat v^2-\n\ell\cd\n v\hat v-\Psi v \hat v+\(\cA\ell_t+\frac{\Psi_t}{2}\)v^2,\\
		\ns\ds \cN\=\ell_t(d\hat v)^2-2\n\ell \cd d(\n v)d\hat v-\Psi dv d\hat v+\ell_t|d(\n v)|^2 +\cA\ell_t(dv)^2+\frac{\Psi_t}{2}(dv)^2,\\
		\ns\ds \cB\=\cA\Psi+(\cA\ell_t)_t-\n\cd \big(\cA \n \ell\big)+\frac{1}{2}\big(\Psi_{tt}-\D  \Psi\big),\\
		\ns\ds V\=2\big(\n\ell\cd\n v\big)\n v-\n\ell |\n v|^2-2\ell_t\n v\hat  v+ \n\ell \hat v^2+\Psi v\n v-\frac{\n\Psi}{2} v^2-\cA v^2\n \ell.
		\ea \right.
		\ee 
		where $(dv)^2$ and $(d\hat v)^2$ denote the quadratic variation processes of $v$ and $\hat v$, respectively.
	\end{lemma}

	Since $T>T^*$, there exists $\k<1$ sufficiently close to $1$ such that $\k T>T^*$. Let
	\bel{as3}
	\ba{ll}
	\ds \a=\frac{\k^2}{n}\min_{x\in \overline G} \frac{\min_{i=1,2\cdots, n}(x_i-x_{0i})^2}{\max_{i=1,2,\cdots,n}(x_i-x_{0i})^2}.
	\ea 
	\ee We have the following lemma.

	\begin{lemma}\label{4.1}
		 There exist constants $c_0, \tilde c_0>0$, $C_0>1$ independent of $ r_2 $
			such that when $ \b\=C_0(1+r_2) $, the following conditions hold: \vspace{1mm}
		
		(1) For any $x\in \overline G$, we have $\ds \max_{i=1,2,\cdots,n}(x_i-x_{0i})^2- \a T^{2}/4<0 $ ;\vspace{1mm}
		
		(2) 
		\bel{w1}
		\ba{ll}
		\min\limits_{i=1,2,\cdots,n}(x_i-x_{0i})^2-n\a^2 \(t-\frac{T}{2}\)^2>c_0>0,\\
		\ns\ds \qq \qq  \forall (t,x)\in \Big\{(t,x)\in Q\|~\sum_{i=1}^n e^{\b(x_i-x_{0i})^2}-ne^{\a\b(t-\frac{T}{2})^2}>0\Big\};
		\ea 
		\ee 

		\vspace{1mm}
		
		(3) 
		
		\begin{equation}\label{11.16-eq2}
			4c_0\b^2+2\b(1-\a)-4r_2\b T>\tilde c_0.
		\end{equation}

	\end{lemma}

	\par {\it Proof of Lemma \ref{4.1}}: First, recalling $\k T>T^*$, and the definition of $T^*$ and $R_1$ given in (\ref{tr1}), we have
	\bel{r1}
	\ba{ll}
	\ds \max_{i=1,2,\cdots,n} (x_i-x_{0i})^2-\a T^2/4\\
	\ns\ds \le \max_{i=1,2,\cdots,n} (x_i-x_{0i})^2-\frac{\k^2T^2}{4n}\min_{x\in \overline G} \frac{\min_{i=1,2\cdots, n}(x_i-x_{0i})^2}{\max_{i=1,2,\cdots,n}(x_i-x_{0i})^2}\\
	\ns\ds < \max_{i=1,2,\cdots,n} (x_i-x_{0i})^2- \frac{T^{*2}}{4n}\min_{x\in \overline G} \frac{\min_{i=1,2\cdots, n}(x_i-x_{0i})^2}{\max_{i=1,2,\cdots,n}(x_i-x_{0i})^2}\\
	\ns\ds \le \max_{i=1,2,\cdots,n} (x_i-x_{0i})^2- R_1^2 \le 0.
	\ea 
	\ee
	Then the assertion (1) follows immediately.

	\par Let
	$$ \tilde c_0=\frac{1-\k^2}{2\k^2}\min_{x\in\overline G}\max_{ i=1,2,\cdots,n}(x_i-x_{0i})^2 $$
	and
	$$c_0=\frac{1-\k^2}{2}\min_{x\in\overline G, i=1,2,\cdots,n}(x_i-x_{0i})^2.$$
	Noting that
	\begin{equation}\label{tr3}
		\lim_{\b\to +\i} \frac{1}{\b}\ln\(\frac{1}{n}\sum_{i=1}^n e^{\b(x_i-x_{0i})^2}\)= \max_{i=1,2,\cdots,n} (x_i-x_{0i})^2
	\end{equation}
	unformly for $x\in\overline G$, we konw that there exists a $\b_0>0$, independent of $x$, such that for any $\b\ge \b_0>1$, it holds that
	\begin{equation}\label{11.16-eq1}
		\frac{1}{\b}\ln\(\frac{1}{n}\sum_{i=1}^n e^{\b(x_i-x_{0i})^2}\)<\max_{i=1,2,\cdots,n} (x_i-x_{0i})^2+\tilde c_0. 
	\end{equation}
	Noting that 
	$$\sum\limits_{i=1}^n e^{\b(x_i-x_{0i})^2}-ne^{\a\b(t-\frac{T}{2})^2}>0,$$
	from \eqref{11.16-eq1}, we know that for $\b\ge \b_0$, 
	\begin{eqnarray}\label{tr2}
		&&\min_{i=1,2,\cdots,n} (x_i-x_{0i})^2-n\a^2\(t-\frac{T}{2}\)^2 \nonumber\\
		&& \ge \min_{i=1,2,\cdots,n} (x_i-x_{0i})^2-n\frac{\a}{\b}\ln \(\frac{1}{n}\sum_{i=1}^n e^{\b(x_i-x_{0i})^2}\)\\
		&& \ge \min_{i=1,2,\cdots,n} (x_i-x_{0i})^2-\k^2\frac{\min_{i=1,2\cdots, n}(x_i-x_{0i})^2}{\max_{i=1,2,\cdots,n}(x_i-x_{0i})^2}\[\max_{i=1,2,\cdots,n} (x_i-x_{0i})^2+\tilde c_0\] \nonumber\\
		&& \ge \frac{1-\k^2}{2}\min_{i=1,2,\cdots,n} (x_i-x_{0i})^2\ge c_0.\nonumber
	\end{eqnarray}
	This implies  that there exists a $\b_0>1$  such that the assertion (2) holds for all $\b\ge \b_0$.
	
	From the definition of the left hand side of \eqref{11.16-eq2}, we know that there exists a constant $C_0>\b_0>1$    such that when $ \b=C_0(1+r_2) $, \eqref{11.16-eq2} holds. We finish the proof of Lemma \ref{4.1}. \endpf 
	
	\ss	
	
	For parameters $\l,\mu>0$ and $ \b $ satisfying conditions in Lemma \ref{4.1}, we choose the weight function $\th$ as follows:
	\bel{weight}
	\ba{ll}
	\ds \th=e^{\ell},\q \ell=\l\phi,\q\phi=e^{\mu\si},
	\ea 
	\ee
	where
	\begin{equation}\label{c3}
		\si(t,x)\= \sum_{i=1}^n e^{\b(x_i-x_{i0})^2}-ne^{\a\b (t-\frac{T}{2})^2}.
	\end{equation}

	Let
	\begin{equation}\label{c4}
		\ba{ll}
		\ds Q(b)\=\big\{ (t,x)\in Q\big| \si(t,x)>b\big\},\q \mbox{ for any constant $b>0$,} \\
		\ns\ds c_1 \=\min_{x\in \overline G} \sum_{i=1}^n e^{(x_i-x_{0i})^2}-n.
		\ea 
	\end{equation}
	Noting that $\b>1$, {
		and  
		$$
		\si\(\frac{T}{2},x\)=\sum_{i=1}^n e^{\b (x_i-x_0)^2}-n>c_1,
		$$
	} we can choose some small constants $\e, \d>0$ with $ \e, \d\sim O(e^{-\b}) $ as $ \b\to +\i $, such that
	\bel{as33}
	Q_0\=\(\frac{T}{2}-\e , \frac{T}{2}+\e \)\t G\subset Q(c_1+2\d)\subset Q(c_1+\d) \subset Q(c_1)\subset (\e, T-\e)\t G\= Q_1.
	\ee 

	\smallskip
	
	Define a cut-off function $\chi\in C^\infty_0(Q)$ as follow:
	
	\bel{av4}
	\chi(t,x)\= \left\{
	\ba{ll}
	\ds 1\qq \mbox{ in } Q(c_1+\d),\\
	\ns\ds 0\qq\mbox{ in }Q\backslash Q(c_1).
	\ea 
	\right.
	\ee 
	Put  
	\bel{chi}
	u=\chi z,\q \hat u=\chi_t z+\chi\hat z,
	\ee 
	where $(z,\hat z)$ is the solution to (\ref{bsystem1}) with $\tau=T$  
	and $(z^T,\hat z^T)\in L_{\cF_T}^2(\Om;H_0^1(G)\t L^2(G))$. 
	From (\ref{bsystem1}), we know $(u,\hat u)$ satisfies that
	\bel{eq1}
	\left\{
	\ba{ll}
	\ds du=\hat u dt+\chi ZdW(t)\q& \mbox{ in }Q,\\
	\ns\ds d\hat u-\D u dt=(\chi_{tt}z+2\chi_t\hat z-2\n \chi\cd\n z-\D \chi  z)dt &\\
	\ns\ds \qq\qq\qq +\chi(a_1z+a_2Z-a_3\hat Z)dt
	+(\chi_tZ+\chi\hat Z)dW(t)&\mbox{ in } Q,\\
	\ns\ds u=0&\mbox{ on } \Si,\\
	\ns\ds u(T)=0,\q  \hat u(T)=0& \mbox{ in }G.
	\ea 
	\right.
	\ee

	In what follows, we shall denote by $ C=C(\Omega, G, T) $ and $ \cC= \cC(\Omega, G, T,\b,r_2)$ with $ \cC\sim O(e^{\b^2}) $ as $ \b\to +\i  $ generic positive constants, which may change from line to line. We have the following Carleman estimate for the system (\ref{eq1}).
	\begin{theorem}\label{Car}
		There exist a positive $\mu_0>0$, such for any  $\mu\geq e^{\b \mu_0}$, there exists $\l_0=e^{C\mu}$ so that for any  $\l\geq\l_0$, and $(u,\hat u)$ satisfying (\ref{eq1}), it holds that
		\begin{eqnarray}\label{bb3}
			&&\l^3\mu^4\dbE \int_Q \phi^3\th^2u^2dxdt+\l\mu \dbE \int_Q \phi \th^2\(|\n u|^2+\hat u^2\)dxdt \nonumber\\
			&&\le \cC\Big\{ \dbE \int_Q \Th\phi  \th^2 \big(\l^3\mu^4\phi^2z^2+|\n z|^2+\hat z^2\big)dxdt\\
			&& \qq+ \l\mu^{3/2}\dbE\int_Q (\Th+\chi^2) \phi \th^2\[\l^2\mu^{5/2}\phi^2|a_4z+Z|^2+|\n(a_4 z+Z)|^2+|a_5\hat z+\hat Z|^2\]dxdt\nonumber \\
			&& \qq +\l\mu\dbE \int_{\Si_0} \phi \th^2\|\frac{\pa z}{\pa \nu}\|^2d\G dt\Big\},\nonumber
		\end{eqnarray} 
		where
		\bel{g1}
		\Th\= |\n \chi|^2+\chi_t^2+\sum_{j,k=1}^n|\chi_{x_jx_k}|^2.
		\ee 
	\end{theorem}

	{\it Proof of Theorem \ref{Car}.} We divide the proof into four steps.
	
	\ss
	
	{\bf Step 1.} By (\ref{weight}) and (\ref{c3}), it is easy to see that
	\bel{a1}
	\ba{ll}
	\ds \ell_{x_j}=2\l\mu\b\phi e^{\b(x_j-x_{0j})^2}(x_j-x_{0j}),\q j=1,2,\cdots,n,\\
	\ns\ds \ell_{x_jx_k} = 4\l\mu^2\b^2\phi e^{\b(x_j-x_{0j})^2+\b(x_k-x_{0k})^2}(x_j-x_{0j})(x_k-x_{0k})\\
	\ns\ds \qq\q +4\l\mu\b^2 e^{\b(x_j-x_{0j})^2}(x_j-x_{0j})^2\d_{jk}+2\l\mu\b \phi e^{\b(x_j-x_{0j})^2} \d_{jk},\q j,k=1,2\cdots,n,\\
	\ns\ds \D \ell = 4\sum_{i=1}^n\[ \l\mu^2\b^2\phi e^{2\b(x_i-x_{0i})^2}(x_i-x_{0i})^2+\l\m \b^2 e^{\b(x_i-x_{0i})^2}(x_i-x_{0i})^2 \]\\
	\ns\ds \qq\q +2\sum_{i=1}^n \l\mu\b \phi e^{\b (x_i-x_{0i})^2},
	\ea 
	\ee 
	\bel{a2}
	\ba{ll}
	\ds \ell_t=-2n\l\mu\a\b \phi e^{\a\b(t-\frac{T}{2})^2}\(t-\frac{T}{2}\),\\
	\ns\ds \ell_{tt}=4n^2\l\mu^2\a^2\b^2\phi e^{2\a\b (t-\frac{T}{2})^2}\(t-\frac{T}{2}\)^2-4n\l\mu\a^2\b^2\phi e^{\a\b(t-\frac{T}{2})^2}\(t-\frac{T}{2}\)^2\\
	\ns\ds \qq -2n\l\mu\a\b \phi e^{\a\b(t-\frac{T}{2})^2},
	\ea 
	\ee
	and 
	\bel{za4}
	\ell_{tx_i}=-4n \l\mu\b\a\b^2\phi e^{\b(x_i-x_{0i})^2+\a\b(t-\frac{T}{2})^2}\(t-\frac{T}{2}\)(x_i-x_{0i}),\q i=1,2,\cdots,n.
	\ee 
	Let us apply  Lemma \ref{lem} to $(u, \hat u)$ and choose $\Psi$ as
	\bel{a3}
	\ba{ll}
	\ds \Psi=-\ell_{tt}+\D \ell-4n\l\mu\a^2\b^2\phi e^{\a\b(t-\frac{T}{2})^2}\(t-\frac{T}{2}\)^2-2n\l\mu\a\b\phi e^{\a\b(t-\frac{T}{2})^2}\\
	\ns\ds \qq -4\l\mu\b^2\phi \sum_{i=1}^n e^{\b(x_i-x_{0i})^2}(x_i-x_{0i})^2-2\l\mu\b\phi \sum_{i=1}^n e^{\b(x_i-x_{0i})^2}.
	\ea 
	\ee 
	Keeping in mind that $\chi \equiv 0$ in $Q\backslash Q(c_1)$ and $ c_0 $ is given in (\ref{w1}), from \eqref{a1}--\eqref{a3} and { recalling that $v$, $\hat v$ are given in Lemma \ref{lem}}, { noting that $\hat v\equiv 0$ in $Q\setminus Q(c_1)$ and  for any $(t,x)\in Q(c_1)$, assertion (2) in Lemma \ref{4.1} holds and
		$$
		\sum_{i=1}^n e^{\b(x_i-x_{0i})^2} > ne^{\a\b(t-\frac{T}{2})^2},
		$$}	
	we have that
	\begin{eqnarray}\label{av5}
		&&\big(\ell_{tt}+\D \ell-\Psi\big)\hat v^2 \nonumber\\
		&&=\Big\{8n^2\l\mu^2\a^2\b^2\phi e^{2\a\b(t-\frac{T}{2})^2}\(t-\frac{T}{2}\)^2 \nonumber\\
		&& \q +4\l\mu\b^2\phi \[ \sum_{i=1}^n e^{\b(x_i-x_{0i})^2}(x_i-x_{0i})^2-n\a^2e^{\a\b(t-\frac{T}{2})^2}\(t-\frac{T}{2}\)^2\]\nonumber\\
		&& \q +2\l\mu\b \phi \(\sum_{i=1}^n e^{\b(x_i-x_{0i})^2}-n\a e^{\a\b(t-\frac{T}{2})^2}\)\Big\}\hat v^2 \nonumber\\
		&&\ge \Big\{8n^2\l\mu^2\a^2\b^2\phi e^{2\a\b(t-\frac{T}{2})^2}\(t-\frac{T}{2}\)^2 \\
		&& \q +4\l\mu\b^2\phi \sum_{i=1}^n e^{\b(x_i-x_{0i})^2}\[\min_{i=1,2,\cdots,n} (x_i-x_{0i})^2-\a^2\(t-\frac{T}{2}\)^2\]\nonumber\\
		&& \q  +2\l\mu\b\phi (1-\a)\sum_{i=1}^n e^{\b(x_i-x_{0i})^2}\Big\}\hat v^2 \nonumber\\
		&& \ge \Big[8n^2\l\mu^2\a^2\b^2\phi e^{2\a\b(t-\frac{T}{2})^2}\(t-\frac{T}{2}\)^2
		+ 4c_0\l\mu\b^2\phi\sum_{i=1}^n e^{\b(x_i-x_{0i})^2}\nonumber\\
		&& \q +2\l\mu\b\phi (1-\a)\sum_{i=1}^n e^{\b(x_i-x_{0i})^2}\Big]\hat v^2,\nonumber
	\end{eqnarray}
	\begin{eqnarray}\label{av6}
		&& \big(\ell_{tt}-\D \ell+\Psi\big)|\n v|^2+2\sum_{j,k=1}^n\ell_{x_jx_k}v_{x_j}v_{x_k} \nonumber\\
		&& =8\l\mu^2\a^2\b^2\phi \[ \sum_{i=1}^n e^{\b(x_i-x_{0i})^2}(x_i-x_{0i})v_{x_i}\]^2\nonumber\\
		&& \q +\Big\{4\l\mu\b^2\phi \[ \sum_{i=1}^n e^{\b(x_i-x_{0i})^2}(x_i-x_{0i})^2-n\a^2e^{\a\b(t-\frac{T}{2})^2}\(t-\frac{T}{2}\)^2\]\\
		&& \qq +2\l\mu\b \phi \[ \sum_{i=1}^n e^{\b(x_i-x_{0i})^2}-n\a e^{\a\b(t-\frac{T}{2})^2}\]\Big\}|\n v|^2\nonumber\\
		&& \ge 8\l\mu^2\a^2\b^2\phi \[ \sum_{i=1}^n e^{\b(x_i-x_{0i})^2}(x_i-x_{0i})v_{x_i}\]^2\nonumber\\
		&& \q + \Big[4c_0\l\mu\b^2\phi\sum_{i=1}^n e^{\b(x_i-x_{0i})}+2\l\mu\b\phi (1-\a)\sum_{i=1}^n e^{\b(x_i-x_{0i})^2}\Big]|\n v|^2,\nonumber
	\end{eqnarray}
	and that
	\bel{bv1}
	\ba{ll}
	-4\n\ell_t\cd \n v\hat v\ge &-8n^2\l\mu^2\a^2\b^2\phi e^{2\a\b(t-\frac{T}{2})^2}\(t-\frac{T}{2}\)^2\hat v^2\\
	\ns&\ds  -8\l\mu^2\a^2\b^2\phi \[ \sum_{i=1}^n e^{\b(x_i-x_{0i})^2}(x_i-x_{0i})v_{x_i}\]^2.
	\ea 
	\ee 
	Combining (\ref{av5})--(\ref{bv1}), we conclude that
	\bel{bv2}
	\ba{ll}
	\ds \big(\ell_{tt}+\D \ell-\Psi\big)\hat v^2+\big(\ell_{tt}-\D \ell+\Psi\big) |\n v|^2+\sum_{j,k=1}^n\ell_{x_jx_k}v_{x_j}v_{x_k}-4\n\ell_t\cd \n v\hat v\\
	\ns\ds \ge \big[4c_0\l\mu\b^2\phi+2\l\mu\b\phi (1-\a)\big]\(\sum_{i=1}^n e^{\b(x_i-x_{0i})}\)(|\n v|^2+\hat v^2).
	\ea 
	\ee 
	\par In what follow, for $k_1,k_2\in \dbN$, we denote by $O(\l^{k_1}\mu^{k_2})$ a function of order $\l^{k_1}\mu^{k_2}$ for large $\l$ and $\mu$. Recalling the definition of $\cA$ in (\ref{VMBA}), we have
	\bel{bv3}
	\ba{ll}
	\ds \cA =\ell_t^2-\ell_{tt}-|\n \ell|^2+\D \ell-\Psi\\
	\ns\ds \q =4\l^2\mu^2\b^2\phi^2\[\a^2 n^2e^{2\a\b(t-\frac{T}{2})^2}\(t-\frac{T}{2}\)^2-\sum_{i=1}^n e^{2\b(x_i-x_{0i})^2}(x_i-x_{0i})^2\]+O(\l\mu^2).
	\ea 
	\ee 
	Consequently, 
	\begin{eqnarray}\label{a5}
		\cB &&\3n\3n=\cA\Psi+(\cA\ell_t)_t-\n\cd(\cA\n \ell)+\frac{1}{2}\big(\Psi_{tt}-\D \Psi\big) \nonumber\\
		&& \3n\3n=\cA_t\ell_t-\n \cA \cd\n \ell +O(\l^3\mu^3)\\
		&&  \3n\3n=32\l^3\mu^4\b^4\phi^3\[\sum_{i=1}^n e^{2\b(x_i-x_{0i})^2}(x_i-x_{0i})^2- \a^2n^2e^{2\a\b(t-\frac{T}{2})^2}\(t-\frac{T}{2}\)^2\]^2  \nonumber\\
		&& \; +O(\l^3\mu^3)+O(\l^2\mu^4).\nonumber
	\end{eqnarray}
	If $ (t,x)\in \big\{(t,x)\in Q|~\si(t,x)>c_1\big\} $, then by (2) in Lemma \ref{4.1}, we have
	\begin{eqnarray}\label{zd1}
		&&\sum_{i=1}^n e^{2\b(x_i-x_{0i})^2}(x_i-x_{0i})^2- \a^2n^2e^{2\a\b(t-\frac{T}{2})^2}\(t-\frac{T}{2}\)^2\nonumber\\
		&&\ge \sum_{i=1}^n e^{2\b(x_i-x_{0i})^2}(x_i-x_{0i})^2-n\min_{i=1,2\cdots,n}(x_i-x_{0i})^2e^{2\a\b(t-\frac{T}{2})^2}\nonumber\\
		&& \ge c_0\(\sum_{i=1}^n e^{2\b(x_i-x_{0i})^2}-ne^{2\a\b(t-\frac{T}{2})^2}\)\\
		&& =c_0\[\(\sum_{i=1}^n e^{2\b(x_i-x_{0i})^2}\)^{1/2}+\sqrt ne^{\a\b(t-\frac{T}{2})^2}\]\[\(\sum_{i=1}^n e^{2\b(x_i-x_{0i})^2}\)^{1/2}-\sqrt ne^{\a\b(t-\frac{T}{2})^2}\]\nonumber\\
		&& \ge c_0\[\(\sum_{i=1}^n e^{2\b(x_i-x_{0i})^2}\)^{1/2}+\sqrt ne^{\a\b(t-\frac{T}{2})^2}\]\( \frac{1}{\sqrt n}\sum_{i=1}^n e^{\b(x_i-x_{0i})^2}-\sqrt ne^{\a\b(t-\frac{T}{2})^2}\) \nonumber\\
		&& \ge \frac{1}{\sqrt n}c_0c_1\[\(\sum_{i=1}^n e^{2\b(x_i-x_{0i})^2}\)^{1/2}+\sqrt ne^{\a\b(t-\frac{T}{2})^2}\]. \nonumber
	\end{eqnarray}
	This, along with \eqref{a5}, implies that
	\bel{zd2}
	\cB v^2\ge \frac{32c_0^2c_1^2}{n}\l^3\mu^4\b^4\phi^3\(\sum_{i=1}^n e^{2\b(x_i-x_{0i})^2}\)v^2 +\big(O(\l^3\mu^3)+O(\l^2\mu^4)\big)v^2.
	\ee 
	By \eqref{zd2}, we know that there exists $ \mu_0>0 $, such for any $ \mu\geq e^{\b\mu_0} $, there exists $ \l_0=e^{C\mu} $, so that for any $ \l\geq\l_0 $, it holds that
	\begin{equation}\label{zd3}
		\cB v^2\ge \frac{16c_0^2c_1^2}{n}\l^3\mu^4\b^4\phi^3\(\sum_{i=1}^n e^{2\b(x_i-x_{0i})^2}\)v^2.
	\end{equation}
	\ms
	
	{\bf Step 2.} The main difficulties to yield Carleman estimate for stochastic hyperbolic-like operator is the estimate of $\cN$ in  (\ref{VMBA}). We focus on this term now. Noting that $v=\th u$ and $\hat v=\th \hat u+\ell_t v$, from \eqref{eq1}, we obtain that
	\bel{e2}
	\ba{ll}
	\ds dv=(\th\chi\ell_t z+\th\chi_t z+\th\chi\hat z) dt+(-a_4v+\cK)dW(t) 
	\ea 
	\ee
	and that
	\begin{equation}\label{e3}
		\ba{ll}
		d\hat v=\big[\th\ell_t\hat u+\ell_{tt}v+\ell_t(\th\chi\ell_t z+\th\chi_t z+\th\chi\hat z)\\
		\ns\ds \qq \q +\th\D u+\th (\chi_{tt}z+2\chi_t\hat z-2\n \chi\cd\n z-\D \chi  z)+\th \chi(a_1z+a_2Z-a_3\hat Z)\big] dt\\
		\ns\ds \qq +\big(-a_5\hat v +\hat{\cK}\big)dW(t),
		\ea 
	\end{equation}
	where
	\bel{e3-1}
	\cK =\th\chi Z+a_4 v =\th\chi(a_4z+Z),
	\ee
	and
	\bel{e3-2}
	\ba{ll}
	\ds \hat{\cK} & =\th\chi_t Z+\th\chi \hat Z+\ell_t\th\chi Z+a_5\hat v\\
	\ns&\ds =\th\chi_t(a_5 z+ Z)+\th\chi(a_5\hat z+\hat Z)+\th\chi\ell_t(a_5 z+Z).
	\ea 
	\ee
	From \eqref{e2}--\eqref{e3-2}, we get that
	\bel{e1}
	\ba{ll}
	\ds \ell_t (d\hat v)^2-2\n \ell \cd (d\n v)d\hat v+\ell_t|d \n v|^2 =\cJ_1+\cJ_2,
	\ea 
	\ee 
	where 
	\bel{e4-1}
	\cJ_1=\big( a_5^2\ell_{t}\hat v^2-2a_4a_5\n \ell\cd \n v\hat v+a_4^2\ell_t|\n v|^2\big)dt,
	\ee 
	and
	\bel{e4-2}
	\ba{ll}
	\ds \cJ_2=\ell_t(\hat \cK^2-2a_5\hat v\hat \cK)dt-2\n\ell\cd \big[ (-\n a_4 v+\n \cK)(-a_5\hat v+\hat \cK)+a_4\hat \cK\n v\big]dt\\
	\ns\ds \qq +\ell_t\big[ |\n a_4v+\n \cK|^2-2a_4\n v\cd (\n a_4 v+\n \cK)\big]dt.
	\ea 
	\ee 
	From \eqref{e4-1}, we see that
	\bel{e6}
	\ba{ll}
	\ds \dbE\int_Q \cJ_1 dxdt &\ds \ge -\dbE\int_Q 2r_2\l\mu\b\phi T \(\sum_{i=1}^n e^{\b(x_i-x_{0i})^2}+ne^{\a\b(t-\frac{T}{2})^2}\)\big(\hat v^2+|\n v|^2\big)dxdt\\
	\ns&\ds\ge -\dbE\int_Q 4r_2\l\mu\b\phi T\( \sum_{i=1}^n e^{\b(x_i-x_{0i})^2}\)\big(\hat v^2+|\n v|^2\big)dxdt.
	\ea 
	\ee
	From \eqref{e4-2}, noting that $ \phi>1 , \forall (t,x)\in Q(c_1)$, we obtain that
	\bel{zd3*}
	\ba{ll}
	\ds \dbE\int_Q \cJ_2dx dt \ge&\ds -\cC \l\mu^{3/2}\dbE \int_Q \phi\big(\hat \cK^2+\phi^2 \cK^2+|\n \cK|^2\big)dxdt\\
	\ns&\ds  -\cC \l\mu^{1/2}\dbE\int_Q \phi \big(\mu\phi^2v^2+\hat v^2+|\n v|^2\big)dxdt.
	\ea 
	\ee 

	Next, we deal with $\ds-\Psi dvd\hat v+\cA \ell_t (dv)^2+\frac{\Psi_t}{2} (dv)^2$.  It follows from (\ref{VMBA}) and (\ref{a3}) that
	\bel{b1}
	\ba{ll}
	\ds \dbE \int_Q \[\Psi dvd\hat v+\cA \ell_t(d v)^2+\frac{\Psi_t}{2}(dv)^2\]dx\\
	\ns\ds =\dbE \int_Q \Psi \big(-a_4v+\cK\big)\big(-a_5\hat v +\hat \cK\big)dxdt+\dbE\int_Q \(\cA+\frac{\Psi_t}{2}\) \ell_t \big(-a_4 v+\cK\big)^2dxdt\\
	\ns\ds \ge -\cC (\l^3\mu^3+\l^2\mu^4)\dbE \int_Q \phi^3(\cK^2+v^2)dxdt-\cC \dbE\int_Q \phi (\hat \cK^2+\hat v^2)dxdt.
	\ea
	\ee
	By \eqref{e3-1} and \eqref{e3-2}, we obtain that
	\bel{abc7}
	\ba{ll}\ds 
	\l\mu^{3/2}\dbE \int_Q \phi \big(\hat \cK^2+\l^2\mu^{5/2} \phi^2 \cK^2+ |\n \cK|^2\big)dxdt\\
	\ns\ds \le \cC\l\mu^{3/2}\dbE \int_Q (\Th +\chi^2)\th^2\phi \big[ \l^2\mu^{5/2}\phi^2(a_4z+Z)^2+|\n (a_4z+Z)|^2+(a_5\hat z+\hat Z)\big]dxdt\\
	\ns\ds \q +\cC\dbE\int_Q (\Th+\chi^2) \l^3\mu^{7/2}\phi^3\th^2 z^2dxdt.
	\ea 
	\ee 

	Combining (\ref{e1}) and \eqref{e6}--(\ref{abc7}), we conclude that
	\begin{eqnarray}\label{e7}
		&& \dbE \int_Q \cN dx\nonumber\\
		&& \ge  -\dbE\int_Q 4r_2\l\mu\b\phi T\( \sum_{i=1}^n e^{\b(x_i-x_{0i})^2}\)(\hat v^2+|\n v|^2)dxdt\nonumber\\
		&&  \q -\cC \dbE \int_Q \big((\l^3\mu^{7/2}+\l^2\mu^4)\phi^3v^2+\l\mu^{1/2}\phi \hat v^2+\l\mu^{1/2}\phi |\n v|^2\big)dxdt \\
		&&  \q-\cC\dbE \int_Q \Th \l^3\mu^{7/2}\phi^3 \th^2z^2dxdt \nonumber \\
		&& \q -\cC\l\mu^{3/2}\dbE \int_Q \big(\Th +\chi^2\big)\th^2\phi \big[ \l^2\mu^{5/2}\phi^2(a_4z+Z)^2+| \n (a_4z+Z)|^2+(a_5\hat z+\hat Z)\big]dxdt.\nonumber
	\end{eqnarray}

	{\bf Step 3.} In this step, we handle the boundary term $d\cM$ and $V$.  Recalling (1) of Lemma \ref{4.1}, we know that the cut-off function $\chi$ satisfies that $\chi(t,x) \equiv 0$ for all $ (t,x)\in Q\backslash Q(c_1)$. This, together with  \eqref{VMBA} and (\ref{as33}), implies 
	\bel{f2}
	\dbE\int_Q dM=0.
	\ee 
	Recalling the definition of $ \Si_0 $ in (\ref{d10}) and $ v=\tilde v=0 $ on $ (0,T)\t \Si $, we get that $ \n v/|\n v|=\nu(x) $ ($ \nu(x)$ represents the unit outward normal vector of $ G $). By (\ref{pointwise_identity}) and (\ref{VMBA}), we have
	\begin{equation}\label{f1}
		\ba{ll}
		\ds \dbE\int_Q \div Vdxdt=\dbE  \int_{\Si} V\cd \nu(x)d\G dt\\
		\ns\ds =  \dbE \int_\Si 2(\n\ell\cd\n v)(\n v\cd \nu)-(\n\ell\cd \nu)|\n v|^2d\G dt\\
		\ns\ds =\dbE \int_\Si (\n\ell\cd\n v)(\n v\cd \nu)d\G dt\\
		\ns\ds \le \cC\l\mu\dbE \int_{\Si_0} \phi \th^2 \|\frac{\pa z}{\pa \nu}\|^2dxdt.
		\ea
	\end{equation}

	{\bf Step 4.} Noting that $(u,\hat u, Z, \hat Z)$ satisfies (\ref{eq1}) and recalling  (\ref{g1}) for $\Th$, we deduce that
	\bel{f5}
	\ba{ll}
	\ds \dbE \int_Q \th\big(-2\ell_t \hat v+2\n\ell\cd \n v+\Psi v\big)\big(d\hat u-\D udt\big)dx \\
	\ns\ds =\dbE \int_Q \th \big(-2\ell_t \hat v+2\n\ell\cd\n v+\Psi v\big)\\
	\ns\ds \qq\qq\times\big[\big(\chi_{tt}z+2\chi_t \hat z-2\n\chi\cd\n z-\D \chi z\big)+\chi\big(a_1z+a_2Z-a_3\hat Z\big)\big]dxdt\\
	\ns\ds \le \frac{1}{2}\dbE\int_Q \th^2 \big( -2\ell_t\hat v+2\n\ell\cd\n v+\Psi v\big)^2dxdt\\
	\ns\ds \q +\cC\dbE \int_Q \Th \th^2 \big(z^2+\hat z^2+|\n z|^2\big)dxdt+\cC  \dbE \int_Q\big(\phi^3v^2+\phi \hat v^2\big)dxdt\\
	\ns\ds \q +\cC \dbE\int_Q \th^2\chi^2 \big( |a_4+Z|^2+|a_5\hat z+\hat Z|^2\big)dxdt.
	\ea 
	\ee 

	At last, recalling that $v=\th \chi z$ and $\hat v=\th \hat u+\ell_t v=\th(\chi_t z+\chi \hat z)+\ell_t \th  \chi z$,  combining  (\ref{11.16-eq2}), (\ref{bv2}), (\ref{zd3}),   (\ref{e7})--(\ref{f5}) together,   we get   (\ref{bb3}).\endpf
	
	\section{Proof of Theorem \ref{1.1}}\label{sec-ob}
	
	This section is devoted to the proof of the observability estimate of \eqref{bsystem1}. We need the following energy estimate.
	\begin{lemma}\label{le}
		For any solution $ (z,\hat z, Z,\hat Z) $ to the equation (\ref{bsystem1}), we have the following inequalities:
		\bel{b5}
		\ba{ll}
		|(z^T,\hat z^T)|_{L_{\cF_T}^2(H_0^1(G)\t L^2(G))}\\
		\ns\ds \le C e^{Cr_2} \dbE \Big\{\int_{Q_0}( \hat z^2+ |\n z|^2+z^2)dxdt \\
		\ns\ds\qq\qq\q  +\int_Q \big[(a_4z+Z)^2+ |\n(a_4z+Z)|^2+(a_5\hat z+\hat Z)^2\big]dxdt\Big\},
		\ea 
		\ee 
		and that
		\bel{b6}
		\ba{ll}
		\ds |(z^T,\hat z^T)|_{L^2_{\cF_T}(H_0^1(G)\t L^2(G))}\ge Ce^{-Cr_2}\dbE \int_Q(\hat z^2+ |\n z|^2+z^2)dxdt.
		\ea 
		\ee 
	\end{lemma}

	Proof of Lemma  \ref{le} follows from some standard energy estimate. We give it here for the sake of completeness.

	{\it Proof of Lemma \ref{le}.  } For any $ t\in (0,T) $, and some constant $ \a$ to be fixed later, by It\^o's formula, we have
	\begin{eqnarray}\label{n1-1}
		&& d(e^{\a t}\hat z^2) + d(e^{\a t}|\n z|^2) + d(e^{\a t }z^2) \nonumber\\[1mm]
		&& =2e^{\a t}\hat z d\hat z + e^{\a t}(d\hat z)^2 + \a e^{\a t}\hat z^2dt + 2e^{\a t}\n z d\nabla z  + e^{\a t}|d\n z|^2  \\[1mm]
		&& \q  + \a e^{\a t} |\n z|^2dt   + 2e^{\a t} z d z + e^{\a t} (dz)^2 + \a e^{\a t}z^2dt. \nonumber
	\end{eqnarray}
	Integrating \eqref{n1-1} on $(t,s)\times G$ and taking mathematical expectation, using integration by parts, we find that
	\begin{eqnarray}\label{n1-2}
		&&\dbE\int_G e^{\a s}\big( \hat z(s)^2 +  |\n z(s)|^2 + z(s)^2\big)dx - \dbE\int_G e^{\a t}\big( \hat z(t)^2 +  |\n z(t)|^2 + z(t)^2\big)dx\nonumber
		\\
		&&= \dbE\int_G \int_t^s\ \big[e^{\a \tau}(d\hat z)^2+\a e^{\a \tau}\hat z^2d\tau\big] dx +\dbE\int_G \int_t^s \big( e^{\a \tau}|d\n z|^2+\a e^{\a \tau} |\n z|^2d\tau\big)dx \nonumber\\
		&&\q +\dbE\int_G \int_t^s\big[e^{\a \tau} (dz)^2+\a e^{\a \tau}z^2d\tau\big]dx +2\dbE\int_G \int_t^s e^{\a \tau}\hat z(d\hat z-\D zd\tau+zd\tau)dx \\
		&& = \a \dbE \int_G \int_t^s e^{\a \tau}\big(\hat z^2+|\n z|^2+z^2\big)dxd\tau  +\dbE \int_G\int_t^s e^{\a \tau} \big(Z^2+ |\n Z|^2+\hat Z^2\big)dxd\tau \nonumber \\
		&& \q + 2\dbE \int_G \int_t^s  e^{\a \tau}\hat z(a_1z+a_2Z-a_3\hat Z+z)dxd\tau.\nonumber
	\end{eqnarray}
	Let $\a=0$. From \eqref{n1-2}, we see that
	\begin{eqnarray}\label{n1-3}
		&&\dbE\int_G \big( \hat z(s)^2 +  |\n z(s)|^2 + z(s)^2\big)dx \nonumber\\
		&& \leq \dbE\int_G \big( \hat z(t)^2 +  |\n z(t)|^2 + z(t)^2\big)dx
		\\ 
		&&\q  +\dbE \int_G\int_t^s \big[ (a_4z+Z)^2+ |\n (a_4z+Z)|^2+(a_5\hat z+\hat Z)^2\big]dxd\tau \nonumber \\
		&& \q + C(1+r_2)\dbE \int_G \int_t^s \big( \hat z(\tau)^2 +  |\n z(\tau)|^2 + z(\tau)^2\big) dxd\tau.\nonumber
	\end{eqnarray}
	It follows from \eqref{n1-3} and Gronwall's inequality that
	\begin{eqnarray}\label{n2}
		&&|(z^T, \hat z^T)|^2_{L^2_{\cF_T(\Om;H_0^1(G)\t L^2(G))}} \nonumber\\
		&& \le e^{C(r_2+1)T}\Big\{\dbE \int_G \(z(t)^2+|\n z(t)|^2+\hat z(t)^2\)dx \\
		&&\qq +C\dbE \int_G\int_t^T \big[ (a_4z+Z)^2+ |\n (a_4z+Z)|^2+(a_5\hat z+\hat Z)^2\big]dxd\tau \Big\},\nonumber
	\end{eqnarray}
	which implies (\ref{b5}).
	
	\ss

	Next,  let $s=T$ in \eqref{n1-2} and $\a=C(r_2+1)$ for a constant $C>0$ such that
	\begin{eqnarray}\label{n1-4}
		&&\Big|2\dbE \int_G \int_t^s  e^{\a \tau}\hat z(a_1z+a_2Z-a_3\hat Z+z)dxd\tau\Big| \nonumber\\
		&&\leq \frac{\a}{2} \dbE \int_G \int_t^s (\hat z^2+|\n z|^2+z^2)dxd\tau\\
		&&\q +\frac{1}{2}\dbE \int_G\int_t^s e^{\a \tau} \big[ Z^2+ |\n Z|^2+\hat Z^2\big]dxd\tau.\nonumber
	\end{eqnarray}
	By \eqref{n1-2} and \eqref{n1-4},  we get that
	\begin{equation}\label{n4}
		\ba{ll}
		\ds |(z^T,\hat z^T)|^2_{L^2_{\cF_T}(\Om;H_0^1(G)\t L^2(G))}\ge  e^{-C(r_2+1)T}\dbE \int_G \(z(t)^2+|\n z(t)|^2+\hat z(t)^2\)dx,
		\ea 
	\end{equation}
	which implies (\ref{b6}).\endpf

	\ms
	
	Now we are in a position to prove the observability estimate.
	
	\ss
	
	{\it Proof of Theorem \ref{1.1}.} By (\ref{as33}) and (\ref{bb3}),  we conclude that
	\bel{b4}
	\ba{ll}
	\ds \dbE\int_{Q_0}(\hat z^2+|\n z|^2+z^2)dxdt\\
	\ns\ds \le \cC \big[\exp(2\l e^{(c_1+\d)\mu_0}) -\exp(2\l e^{(c_1+2\d)\mu_0})\big]\dbE \int_Q \big(\hat z^2+|\n z|^2+z^2\big)dxdt\\
	\ns\ds \q  + \cC\exp(2\l e^{\mu_0e^{C\b}})\Big\{\dbE \int_{\Si_0}\|\frac{\pa z}{\pa \nu}\|^2d\G dt\\
	\ns\ds \qq\qq\qq\qq\q   +\dbE \int_Q \big[ (a_4z+Z)^2+ |\n (a_4z+Z)|^2+(a_5\hat z+\hat Z)^2\big]dxdt \Big\}.
	\ea 
	\ee

	\par 
	
	It follows from (\ref{b5}), \eqref{b6} and \eqref{b4} that 
	\begin{eqnarray}\label{b7}
		&&|(z^T,\hat z^T)|_{L^2_{\cF_T}(\Om;H_0^1(G)\t L^2(G))}\nonumber\\
		&& \le \cC e^{Cr_2}\[\exp\big(2\l e^{(c_1+\d)\mu_0}\big) -\exp\big(2\l e^{(c_1+2\d)\mu_0}\big)\]|(z^T,\hat z^T)|_{L^2_{\cF_T}(\Om;H_0^1(G)\t L^2(G))}\nonumber\\
		&&\q +\cC\exp\big(2\l e^{\mu_0e^{C\b}}+Cr_2\big)\Big\{\dbE \int_{\Si_0}\|\frac{\pa z}{\pa \nu}\|^2d\G dt\\
		&& \qq\qq +\dbE \int_Q \[ (a_4z+Z)^2+ |\n(a_4z+Z)|^2+(a_5\hat z+\hat Z)^2\]dxdt\Big\}.\nonumber
	\end{eqnarray} 
	Let us choose $\l$ large enough such that
	$$\cC e^{Cr_2}\[\exp(2\l e^{(c_1+\d)\mu_0}) -\exp(2\l e^{(c_1+2\d)\mu_0})\]<1.$$
	Recalling  Lemma \ref{4.1}  that $ \b=C_0(1+r_2) $  and notice that $ \d  \sim O(e^{-\b})$ as $ \b \to +\i $, we obtain the inequality (\ref{3}). \endpf


\begin{thebibliography}{}
		
		
		\bibitem{AI} S.~A.~Avdonin and S.~A.~Ivanov. \it
		Families of Exponentials. The Method of Moments
		in Controllability Problems for Distributed
		Parameter Systems. \sl  Cambridge University
		Press, Cambridge. \rm 1995.
		
		\bibitem{barbu1} V.~Barbu, A.~R$\breve{\rm a}$\c scanu and G.~Tessitore.
		\it Carleman estimate and cotrollability of
		linear stochastic heat equatons. \sl Appl. Math.
		Optim. \rm {\bf 47} (2003), 97--120.
		
		\bibitem{Bardos-Lebeau-Rauch1} C.~Bardos, G.~Lebeau and J.~Rauch.
		\it Sharp sufficient conditions for the
		observation, control and stabilizion of waves
		from the boundary. \sl SIAM J. Control Optim.
		\rm {\bf 30} (1992), 1024--1065.
		
		
		\bibitem{JMC} J. -M. Coron. \it Control and Nonlinearity. \sl  American Mathematical Society, Providence, RI. \rm  2007.
		
		\bibitem{minicourse} R. Dalang, D. Khoshnevisan, C. Mueller, D. Nualart and Y. Xiao. \it A Minicourse on Stochastic Partial Differential Equations. \sl Springer-Verlag, Berlin. \rm  2009.
		
		\bibitem{FL1}  X. Fu and X. Liu. \it  Controllability and observability of some stochastic complex Ginzburg-Landau equations. \sl SIAM J. Control Optim. \rm {\bf 55} (2017),  1102--1127. 
		
		
		
		\bibitem{GCL} P.~Gao, M.~Chen and Y.~Li. \it Observability
		estimates and null controllability for forward
		and backward linear stochastic
		Kuramoto-Sivashinsky equations. \sl SIAM J.
		Control Optim. \rm{\bf 53} (2015), 475--500.
		
		
		\bibitem{JLL} J. -L. Lions. \it  Contr\^olabilit\'o  Exacte, Perturbations et Stabilisation de Syst\`emes Distribu\'es.  Tome 1, Contr\^olabilit\'e Exacte. \sl Masson, Paris. \rm 1988.
		
		\bibitem{LiuY} X.~Liu and Y. Yu. \it   Carleman estimates of some stochastic degenerate
		parabolic equations and application. \sl SIAM J. Control Optim.
		\rm{\bf 57} (2019),  3527--3552.
		
		\bibitem{Lu0}  Q.~L\"{u}. \it Some results on the
		controllability of forward stochastic heat
		equations with control on the drift. \sl J.
		Funct. Anal. \rm {\bf 260} (2011), 832--851.
		
		
		
		\bibitem{Lu4} Q. L\"{u}. \it Exact controllability for stochastic Schr\"odinger
		equations. \sl J. Differential Equations.
		\rm{\bf 255} (2013), 2484--2504.
		
		\bibitem{Luqi8} Q. L\"{u}. \it  Exact controllability for stochastic
		transport equations. \sl SIAM J. Control Optim.
		\rm{\bf 52} (2014),  397--419.
		
		\bibitem{Lu2022b} Q. L\"{u}. \it  Control theory of stochastic distributed parameter
		systems: recent progresses and open problems. \sl  Proceedings of the International Congress of Mathematicians. \rm   2022. In press.
		
		\bibitem{Lu2022} Q. L\"{u} and Y. Wang. \it Null controllability for fourth order stochastic parabolic equations. \sl SIAM J. Control Optim. \rm {\bf 60} (2022),  1563--1590.
		
		\bibitem{QLXZ1} Q. L\"u and X. Zhang. \it Exact controllability for a refined stochastic wave equation.  \rm  airXiv: 1901.06074.
		
		\bibitem{QLXZ} Q. L\"u  and  X. Zhang. \sl Mathematical Control Theory for Stochastic Partial Differential Equations. \rm Springer-Verlag. \rm   2021.
		
		\bibitem{EN} E. Nelson. \it Dynamical Theories of Brownian Motion. \sl Princeton University Press, Princeton, N.J.  \rm 1967.
		
		
		\bibitem{DLR} D. L. Russell. \it Controllability and stabilizability theory for linear partial differential equations: recent progress and open problem. \sl SIAM Rev. \rm {\bf 20} (1978), 639--739.
		
		\bibitem{TZ} S.~Tang and X.~Zhang. \it Null controllability
		for forward and backward stochastic parabolic
		equations. \sl SIAM J. Control Optim. \rm{\bf 48} (2009), 2191--2216.
		
		\bibitem{Yan2022} L. Yan,  B. Wu,  S. Lu and Y. Wang. \it Null controllability and inverse source problem for stochastic Grushin equation with boundary degeneracy and singularity. \sl ESAIM Control Optim. Calc. Var. \rm{\bf 28} (2022), Paper No. 43, 34 pp.
		
		\bibitem{Yu2022} Y. Yu and J. Zhang. \it Carleman estimates of refined stochastic beam equations and applications. \sl SIAM J. Control Optim. \rm{\bf  60} (2022),  2947--2970. 
		
		\bibitem{Zhangxu4} X.~Zhang.  \it Explicit observability inequalities
		for the wave equation with lower order terms by
		means of Carleman inequalities. \sl SIAM J.
		Control Optim.  \rm {\bf 39}  (2000),  812--834.
		
		\bibitem{XZ} X. Zhang. \it  A unified controllability/observability theory for some stochastic and deterministic partial differential
		equations. \sl  Proceedings of the International Congress of Mathematicians.\rm  Volume IV, 3008--3034, Hindustan Book
		Agency, New Delhi, 2010.
		
		\bibitem{EZ} E. Zuazua. \it Controllability and observability of partial differential equations: some result and open problems. \sl in Handbook of Differential Equations: Evolutionary Equations. \rm {\bf 3} (2006), 527--621
		
	\end{thebibliography}
\end{document}